\documentclass[english, 12pt]{article}
\usepackage[english]{babel,varioref}
\usepackage[ansinew]{inputenc}
\usepackage[fleqn,reqno]{amsmath}
\usepackage{amsthm}
\usepackage{paralist}
\usepackage{natbib}
\usepackage{graphics}
\usepackage{epic,eepic}
\usepackage[doublespacing]{setspace}
\usepackage{lineno}

\bibpunct[: ]{(}{)}{,}{a}{}{,}

 \setdefaultenum{(i)}{a)}{(1)}{(A)}

\theoremstyle{plain}

\newtheorem{thm}{Theorem}
\newtheorem{cor}[thm]{Corollary}
\newtheorem{lem}[thm]{Lemma}

\newtheorem{defn}[thm]{Definition}

\def\ld{\mathop{\rm ld}}
\begin{document}


\title{Toward a deeper understanding of a basic cascade}

\author{Uwe Saint-Mont,
Nordhausen University of Applied Sciences}

\maketitle

\pagebreak
\section*{Abstract}

Towards the end of the last century, B. Mandelbrot saw the importance, revealed the beauty, and robustly promoted (multi-)fractals. Multiplicative cascades are closely related and provide simple models for the study of turbulence and chaos.

For pedagogical reasons, but also due to technical difficulties, continuous stochastic models have been favoured over discrete cascades.
Particularly important are the $\alpha$ and the $p$ model.
It is the aim of this contribution to introduce original concepts that shed new light on a variant of the latter paradigmatic cascade and allow key features to be derived in a rather elementary fashion.

To this end, we introduce and study a discrete version of the $p$ model which is based on a new kind of sampling. Technical machinery can be kept simple, therefore proofs are straightforward and formulas are explicit. It is hoped that the proposed line of investigation may enhance understanding and simplify received multifractal analyses.

\vspace{3ex}
{\bf Keywords}: cascades, $p$ model, multifractals, sampling, law of large numbers, Mandelbrot, Bernoulli


\vspace{5ex}
\section*{Acknowledgements}

In 2019, I submitted an earlier version of this manuscript to {\it Nonlinear Processes in Geophysics}. The author would like to thank the editor of that journal, Daniel Schertzer, and several anonymous referees for very valuable suggestions that improved the paper considerably.

\newpage

\begin{doublespace}

\section{Introduction}\label{applied}

Early examples of fractals were provided by, among others, mathematicians Weierstra\ss, Cantor and Peano. Later, upon studying dynamic systems, chaos and turbulence, physicists found similar patterns. Pioneering work dates back to the first half of the 20th century, in particular to \citet{ri22} and \citet{ko41}, and \citet{ma82,ma97,ma99} obtained a first synthesis when he established a strong link between fractal geometry and its applications in the sciences (physics and economics in particular) which has since been extended to ``multifractal methodology'' \citep{sa17}, general ``critical phenomena'' \citep{so06}, and asymptotic theory \citep{ke11}.

Moving from the objects involved to the processes generating them, cascades have come into focus \citep{sc11} only recently. On the one hand, they are quite common. On the other hand,
they are also a ``key idea'' conceptually (\citet{lo19}, p. 76). That is, although cascades are rather primitive (iterate a basic building block), they can easily be adapted to observable phenomena: The basic building block (a type of fork-like structure) can be chosen appropriately, the propagation mechanism may be deterministic or stochastic, scales (and scale invariance) are closely related, and the approach works in spaces of (almost) any dimension.

Endowed with an abundance of nonlinear processes, geophysics has been especially productive in this respect (see, in particular, \citet{ma89}, and \citet{lo07}), with contributions ranging from the atmosphere (climate and weather), wave dispersion and topography to geology and mining (e.g., \citet{se10}, \citet{lo13}, \citet{ag19}).

In more detail, the basic element to be considered is a bifurcation which distributes some commodity $m$ (mass, energy,$\ldots$) in a parent vertex among two descendants:
$$
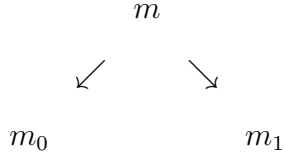
\begin{array}{ccccc}
 && m &&  \\
 &\swarrow &  & \searrow & \\
 m_0 &&&& m_1 \\
\end{array}
$$
 If mass is neither added nor withdrawn (i.e., conserved, $m_0+m_1=m$), such a step is called {\it microcanonical}. Typically, there is a preference toward one of the siblings. Thus \citet{de51, de53} defined $m= 2x, m_0= (1-d)x$ and $m_1=(1+d)x$, where $0 \le d \le 1$. (For a recent application, see, e.g., \citet{ag19}.). \citet{sa17}, p. 471, start with the uniform distribution on the unit interval $U(0,1)$, and allocate the mass in $(1/3,2/3)$ to the interval $(2/3,1)$. Thus $m=1$, the mass in the left-hand interval (0,1/3) remains $m_0=1/3$, and the mass in the right-hand interval $(2/3,1)$ increases to $m_1=2/3$.

 The standard stochastic model for a bifurcation (the division of mass among two descendants) is the Bernoulli distribution $B(p)$, with $m=1, m_0=1-p$ and $m_1=p$. In the following, we are going to study the more general `double random' model
$$
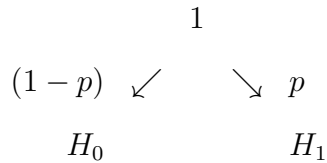
\begin{array}{rcccl}
 && 1 &&  \\
(1-p) &\swarrow &  & \searrow & p \\
 H_0 &&&& H_1 \\
\end{array}
$$
That is, a Bernoulli random variable $C \sim B(p)$ selects population $V_1 \sim H_1$ with probability $p$ and population $V_0 \sim H_0$ with probability $1-p$. Of course, if $V_i \equiv i$ $(i=0,1)$ we just have a standard Bernoulli with realizations 0 and 1, respectively. However, in general, $V_0$ and $V_1$ are random variables with probability distributions $H_0$ and $H_1$, respectively. Thus the bifurcation turns out to be canonical in the sense that mass is conserved {\it in expectation}: If, w.l.o.g., $EV_i = i$, then $p EV_1 + (1-p) EV_0 = p.$ However, $p v_1 + (1-p) v_0$ need not be $p$.

The second constitutive idea of a cascade is an iteration of the splitting procedure. That is, every leaf of a given structure (a tree, if there is a single root) is replaced by another bifurcation. Although we are used to think of the Binomial $B(n,p)$ as the most natural continuation of the Bernoulli $B(p)$, note that there is a fundamental difference between a cascade and a binomial structure (see the next figure).

\vspace{1ex}
{\bf Illustration 1:}\label{illu1} Cascade (local splitting) vs. binomial structure.
$$
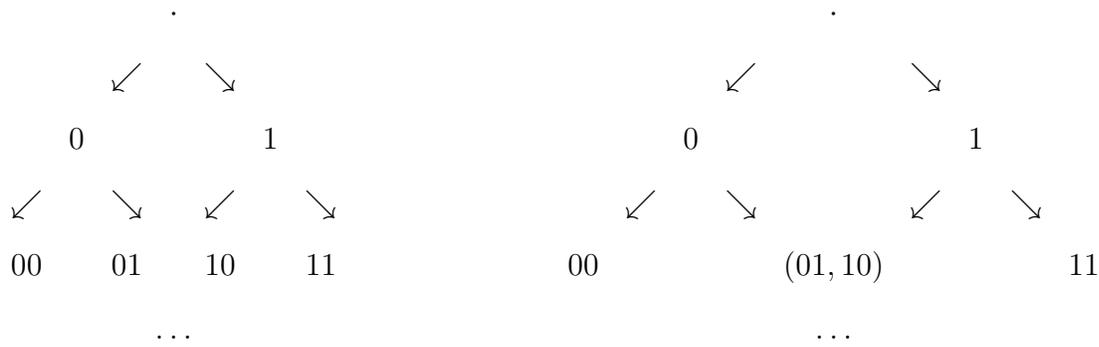
\begin{array}{lcccccr}
&&& \cdot &&& \\
 &  & \swarrow && \searrow &  & \\
 & 0 &&&& 1 & \\
  \swarrow &  & \searrow &  & \swarrow & & \searrow \\
  00 &  & 01 &  & 10 &  & 11 \\
  \multicolumn{7}{c}{\ldots}
\end{array}
\hspace{15ex}
\begin{array}{lcccccccr}
&&&& \cdot &&&& \\
 &  && \swarrow && \searrow & & & \\
 && 0 &&&& 1 && \\
 & \swarrow &  & \searrow &  & \swarrow & & \searrow & \\
  00 & & &  & (01, 10) & & &  & 11 \\
  \multicolumn{9}{c}{\ldots}
\end{array}
$$
In a cascade, there are only bifurcations. Therefore, a binary tree evolves whose number of leafs doubles with every iteration and grows exponentially fast $(2^n \rightarrow 2^{n+1}, n \ge 0)$. A binomial structure, however, grows slowly $(n \rightarrow n+1, n\ge 1)$, since right after the splits, descendants merge. Actually, given $n$ iterations, the characteristic feature of the Binomial is to count the number of ways $\binom{n}{k}$ that lead to some leaf $k$, $(0 \le k \le n)$.

 From a dynamic system perspective, paths split {\it and} merge in Pascal's triangle, leading to an accumulation of mass toward the centre. However, since paths only split in a cascade, the latter are excellent models for divergent phenomena, such as chaos and turbulence (with vortices, eddies or boxes multiplying, but not fusing). Moreover, given a starting point with all mass concentrated there, repeated local bifurcations (all governed by the same mechanism), evoke self-similar structures that can often be extended to a reasonable (multi-)fractal limit \citep{sh19}.

It is also typical that the division of the available mass is not `fair'. For instance, the cascade considered by \citet{sa17}, p. 471, favours the right-hand side with factor $f = 2:1$. Thus their figure at the top of that page can be translated into the numbers
$$
\begin{array}{cccccccccccccccc}
   &  &  &  & &&& 27  &&&& & & &  &  \\\hline
   & & &  9 & && & & &&& 18  & &  &  \\\hline
   & 3 & &  && 6 & &  && 6 &  & && 12 &  \\\hline
  1 && 2 && 2 && 4  &  & 2 && 4 && 4 && 8 \\\hline
\end{array}
$$

Trained to detect the Binomial, one may readily see $\binom{3}{k}$ times the number $2^k$ $(0 \le k \le 3)$ in the last line. However, that is just a summary.
More importantly, there is a self-similar structure in every line, in particular $2 \cdot 1 = 2$; $2 \cdot (1,2) = (2,4)$ and $ 2 \cdot (1,2,2,4) = (2,4,4,8)$ in the last.

Pascal's arithmetic triangle corresponds to the Binomial $B(n,p)$. Analogously, the above `geometric triangle' (with $f=2/1=p/(1-p)$, and thus $p=2/3$) corresponds to a new probability distribution that will be named the Weaver's distribution $W(n,p)$. More explicitly, $Y \sim W(3,2/3)$ and $W(3,p)$, respectively, obtain the values $y_k$ with probabilities $p_k$:
$$
\begin{array}{|c|c|c|c|c|c|c|c|c|}\hline
y_k & 0 & 1/7 & 2/7 & 3/7 & 4/7 & 5/7 & 6/7 & 1 \\\hline
p_k & 1/27 & 2/27 & 2/27 & 4/27 & 2/27 &  4/27 &  4/27 & 8/27 \\
p_k &  (1-p)^3 & p(1-p)^2 & p (1-p)^2 & p^2 (1-p) & p (1-p)^2 &  p^2 (1-p) &  p^2 (1-p) & p^3 \\ \hline
\end{array}
$$

The choice of the realizations $y_k = k/(2^n-1)$ or the unstandardized values $k=0,1,\ldots,2^n-1$, is a very natural one, and simplifies a thorough analysis considerably.

Actually, the received $p$ model (see \citet{ma89}, p. 19, \citet{ma74}, p. 329, and \citet{de51, de53}) is a continuous version of the Weaver. That is, for $n=0$ one may start with the continuous Uniform distribution on the unit interval.
Next, the proportion $1-p$ is uniformly distributed on the interval $(0,1/2)$, and the proportion $p$ is uniformly distributed on the interval $(1/2,1)$. In the same vein, one splits the masses further (locally), and obtains the following cascade (illustration 2).

\vspace{1ex}
{\bf Illustration 2:}\label{illu-p} Cascade of the $p$ model
$$
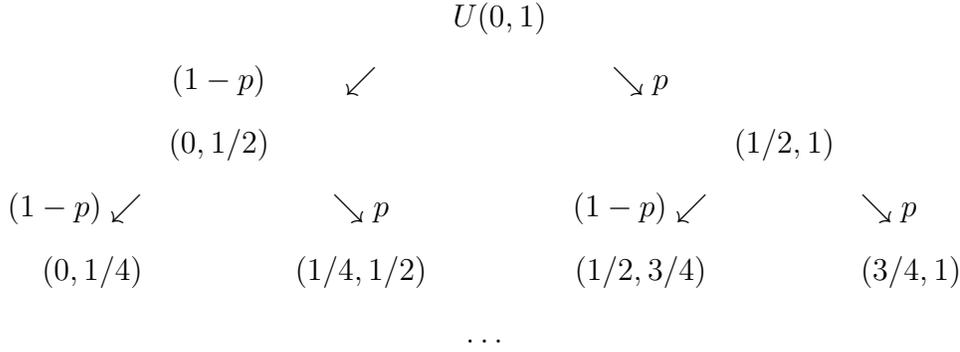
\begin{array}{rcccccl}
&&& U(0,1) &&& \\
 & (1-p) & \swarrow && \searrow p &  & \\
 & (0,1/2) &&&& (1/2,1) & \\
(1-p)  \swarrow &  & \searrow p &  & (1-p) \swarrow & & \searrow p \\
  (0,1/4) &  & (1/4,1/2) &  & (1/2,3/4) &  & (3/4,1) \\
  \multicolumn{7}{c}{\ldots}
\end{array}
$$

For $n=3$, this cascade yields
a piecewise defined density $f(z)$ on the intervals $I_k=(k/8,(k+1)/8)$, where $k=0,\ldots,7$:\footnote{For arbitrary $n$ and $p$ see www.demonstrations.wolfram.com, ``MandelbrotsBinomialMeasureMultifractal.''}
$$
\begin{array}{|c|c|c|c|c|c|c|c|c|}\hline
I_k & \left(0,\frac{1}{8} \right) & \left( \frac{1}{8}, \frac{2}{8} \right) &  \left( \frac{2}{8}, \frac{3}{8} \right)  &  \left( \frac{3}{8}, \frac{4}{8} \right)  &  \left( \frac{4}{8}, \frac{5}{8} \right)  &  \left( \frac{5}{8}, \frac{6}{8} \right)  &  \left( \frac{6}{8}, \frac{7}{8} \right)  &  \left( \frac{7}{8}, 1 \right)  \\\hline
f(z) & 8 (1-p)^3 & 8 p(1-p)^2 & 8 p (1-p)^2 & 8 p^2 (1-p) & 8 p (1-p)^2 & 8 p^2 (1-p) & 8 p^2 (1-p) & 8 p^3 \\ \hline
\end{array}
$$

Curiously enough, \citet{ma99}, p. 87, says that the $p$ model appeared ``in an esoteric corner of mining engineering science.'' However, if one thinks about it, ores are the result of an enrichment process, and a straightforward model for such a process is a sequence of binary decisions, i.e., a cascade that is {\it biased} in favour of some mineral.

Although a density with many jumps is more difficult to handle than a suitably defined distribution with finite support, it is well-known that the limit distribution function of the $p$ model, with the exception of $p=1/2$, has no density \citep{sa43}.
  Instead, one encounters a multifractal structure asymptotically, involving a certain amount of polarization, depending on $p$. This corresponds to veins of gold or a mineral deposit vs. dead rock, say, in a mining application. For a graphic example see \citet{hi99}.

Since the basic building block used in the $p$ model is a binary bifurcation (each point bequeaths its mass to two descendants with proportions $p$ and $1-p$, respectively), the corresponding cascade should be named after {\it Bernoulli}. Unfortunately, the terms `binomial cascade' (and `binomial measure') have caught on in the literature, since a process that is based on the Bernoulli readily yields a binomial distribution.
However, the Weaver $W(n,p)$ and its limit $W(p)$
are as basic as the Binomial and its siblings.

We are going to demonstrate that the crucial difference between the Binomial and the Weaver can be reduced to a slightly more sophisticated way of sampling that `augments' Pascal's triangle to the `geometric triangle' above. The latter multiplicative pattern is equivalent to local Bernoulli bifurcations and brings out the fractal nature of zero-one decisions.
 In algorithmic terms:
 \begin{enumerate}
   \item[0)] There are two populations $H_0,H_1$; fix $n \ge 1$, and let $j=1$.
   \item[1)] Choose one of the populations at random, $C \sim B(p)$.
   \item[2)] Draw an iid sample of size $2^{j-1}$ from the selected population $H_c$.
   \item[3)] Let $j=j+1$.
   \item[4)] If $j \le n$, proceed to step 1)
 \end{enumerate}
Vividly, step 2 (named {\it exponential sampling} hereafter) acts as a ``repulsive force'' that is able to prevent the two distinct populations from merging. Essentially, we are going to study the sum of the sample $X_1;X_2,X_3;X_4,X_5,X_6,X_7;\ldots$ thus constructed.\footnote{Of course, if the populations are trivial, in particular if $H_c \stackrel{d}{=} \varepsilon_c$ (point mass at $c, c \in \{0,1\})$, then the sample just consists of $n$ blocks of zeros and ones, the length of each block being $2^{j-1}$.}

The rest of this article is organized as follows: In the next section we give a rather different, theoretical motivation for the above model. Section 3 defines exponential sampling, and derives the Weaver's distribution and its properties in a systematic way. It turns out that the corresponding deterministic cascade forms `threads' that interweave in a particular way.
Moments are given in section 4, the limit distribution is discussed in section 5, and section 6 is devoted to populations with finite variances. In particular, the variance can be decomposed into several components. Finally, in section 7, we compare several kinds of asymptotic behaviour.

\section{Theoretical motivation}

Given an iid sequence $X_1,X_2,\ldots$ of random variables, the basis of traditional (Frequentist) statistics are the central limit theorem (CLT) and some law of large number (LLN), i.e., the convergence of ${\bar X}_n=S_n /n=\sum_{i=1}^n X_i/n$ towards a single number. However, in calculus, convergence of a sequence $x_1,x_2,\ldots$ is a strong assumption, and, typically, not even the (much weaker) Cesaro-limit $\lim_{n \rightarrow \infty}{\bar x}_n= \lim_{n \rightarrow \infty}(\sum x_i/n)$ exists. In dynamic system theory, also, convergence towards a point is a rare exception.

In probability theory, the iid model represents a single population and a large, potentially infinite sample from this population. To {\it avoid} convergence, it is thus straightforward to consider {\it two} populations (distributions $H_0$ and $H_1$), and a sample that fluctuates between them. In other words, if one switched between the populations skilfully, ${\bar X}_n$ should not converge. In the jargon of dynamic system theory, the (unique) limit point may be replaced by a (more complicated) attractor.

However, a constant switching rate won't do: If $j$ observations from $H_0$ are followed by $j$ observations from $H_1$, and so forth, the arithmetic mean of this sequence will converge, since the `influence' of another $j$ observations on ${\bar X}_n$ becomes insignificant with increasing $n$. Yet if $2^j$ $(j \ge 0)$ observations from $H_0$ are followed by $2^{j+1}$ observations from $H_1$, etc., one then obtains the desired effect. (On a logarithmic scale, taking $\ld =\log_{2}$, the ratio $\ld ((2^{j+1} ) / 2^{j}) = j+1-j =1$ is a constant. Thus, there, one switches at a constant rate, `1' indicating that $H_0$ alternates with $H_1$.) Since $2^0 + 2^2 +2^4+\ldots$ observations are from $H_0$, and $2^1+2^3+2^5+\ldots$ observations are from $H_1$, given a sample of size $2^n-1$, considerably more than one half of these observations come from $H_0 (H_1)$, if $n$ is an odd (even) number. Thus the arithmetic mean cannot `settle' in some point.

Altogether, we obtain a stochastic process that is inhomogeneous in a particular way. Its paths depend on the concrete distributions of $H_0$ and $H_1$, and on the way switching is done. The aim of this article is to explore straightforward consequences of this basic setting.

\section{The weaver's distribution}

In order to keep things finite, suppose for the rest of this contribution that first moments exist, such that without real loss of generality\footnote{See the remark after the proof of Theorem \ref{weaver}, p. \pageref{simplify}} $\mu(H_0)=0$ and $\mu(H_1)=1$ are the expected values of the two populations (distributions) involved.

A particularly simple way to alternate between $H_0$ and $H_1$ is to take the next batch of $2^j$ observations $(j=0,1,\ldots)$ from population $H_0$ with probability $1-p$, and from population $H_1$ with probability $p$. To avoid trivialities, we assume $0<p<1$ throughout this contribution. Thus, one creates a hierarchical random system (a particular random probability measure) composed of a choice mechanism which selects the population in charge, and a realization mechanism which provides observations from the population selected.

\begin{defn}\label{progsamp}
  (Exponential sampling)

  Given two parent distributions $H_0$ and $H_1$, and $0<p<1$, define exponential sampling as follows:

  A sample of size $2^n -1$, i.e., $X_1; X_2, X_3; X_4, X_5, X_6 , X_7;\ldots; X_{2^{n-1}},\ldots,X_{2^n-1}$, consists of $n$ sub-samples, where sub-sample $j$ runs from $X_{2^{j-1}}$ to $X_{2^{j}-1}$ $(j=1,\ldots,n)$.
Further suppose that ${\bf B}_n=(B_{n-1},\ldots,B_0)$ is a vector of $n$ independent Bernoullis $B_j \sim B(p)$.

Then sub-sample $j$ consists of $2^{j-1}$ iid random variables which are distributed according to $H_1$ if $B_{j-1}=1$ and according to $H_0$ if $B_{j-1}=0$.

\end{defn}

With probability $p$, the first observation comes from $H_1$, and with probability $1-p$, the first observation comes from $H_0$. Thus, conditional on this choice, the expected value observed is either $\mu(H_1)=1$ or $\mu(H_0)=0$, and the unconditional mean is $\mu=p \mu(H_1) + (1-p) \mu(H_0) =p$.

With probability $p$, the second {\it and} third observations both come from $H_1$, and with probability $1-p$, these observations both come from $H_0$. Thus, after two choices, the overall situation is as follows:

\begin{tabular}{|l|l|l|l|}
  \hline
   Number of observations & Number of observations & Probability & Conditional \\
   from $H_0$ & from $H_1$ & & Mean \\\hline
    1+2 & 0 & $(1-p)^2$ & 0 \\
    1 & 2 & (1-p) p & $2/3$ \\
    2 & 1 & p (1-p) & $1/3$ \\
    0 & 3 & $p^2$ & 1 \\
  \hline
\end{tabular}

 The unconditional mean does not change, since
$$
\mu = p^2 + \frac{1}{3}p(1-p) + \frac{2}{3} (1-p)p = p^2 +p (1-p) =p.
$$

Similar to the binomial distribution, every path, being defined by a sequence of independent binary decisions $(b_{n-1},\ldots,b_0)$, splits upon moving from $n$ to $n+1$. However, unlike the Binomial, the paths do not merge. Rather, like threads, they interweave (see the next figure):

$$
\begin{array}{rrrrrrrr}
&&& 0 & 1 & && \\
&& {\bf 0} 0 & {\bf 0}  1 & {\bf 1}  0  & {\bf 1}  1 &   &  \\
 \underline{0}  {\bf 0}  0  & \underline{0} {\bf 0} 1  &  \underline{0} {\bf 1} 0 &  \underline{0} {\bf 1} 1 &  \underline{1} {\bf 0}  0 &  \underline{1} {\bf 0} 1 &  \underline{1} {\bf 1} 0 &  \underline{1} {\bf 1}  1  \\
 \multicolumn{8}{c}{\ldots}
\end{array}
$$

The next illustration demonstrates that, in a sense, the difference between weaving and splitting is minor: Given a binary string, weaving adds the next cipher to the left (a prefix), whereas splitting adds the next cipher to the right (a suffix).

\vspace{1ex}
{\bf Illustration 3:}\label{illu2} Global weaving (left) and local splitting (a cascade, right)
$$
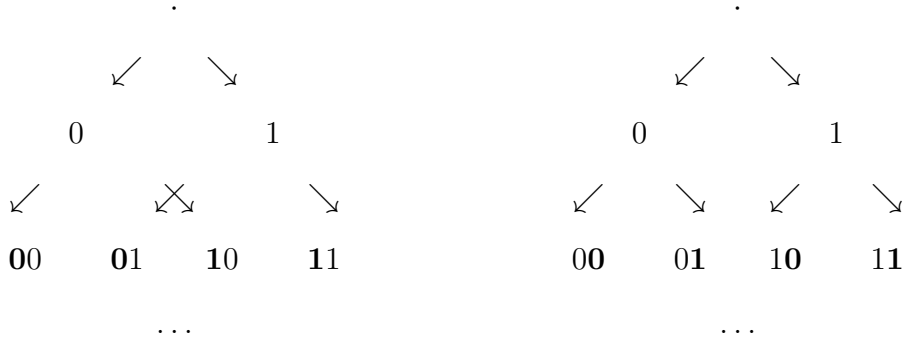
\begin{array}{ccccccc}
&&& \cdot &&& \\
 &  & \swarrow && \searrow &  & \\
 & 0 &&&& 1 & \\
  \swarrow & \multicolumn{5}{c}{\searrow \hspace{-3ex} \swarrow}  & \searrow \\
 {\bf 0}0 &  & {\bf 0} 1 &  & {\bf 1 } 0 &  & {\bf 1} 1 \\
    \multicolumn{7}{c}{\ldots}
\end{array}
\hspace{15ex}
\begin{array}{ccccccc}
&&& \cdot &&& \\
 &  & \swarrow && \searrow &  & \\
 & 0 &&&& 1 & \\
  \swarrow &  & \searrow &  & \swarrow & & \searrow \\
  0 {\bf 0} &  & 0 {\bf 1} &  & 1 {\bf 0} &  & 1 {\bf 1} \\
  \multicolumn{7}{c}{\ldots}
\end{array}
$$

After $n$ steps (selections, choices), one thus obtains an interesting distribution:

\begin{thm}\label{weaver}
              (The weaver's distribution)

 Given the situation described in Definition \ref{progsamp}, suppose the first moments are $\mu(H_0)=0$ and $\mu(H_1)=1$, respectively.  

For $n=1,2,\ldots$ let $S_n=\sum_{i=1}^{2^n-1} X_i$, ${\bar X}_n=S_n / (2^n-1)$, and $Y_n = E({\bar X}_n | {\bf B}_n)$. Some elementary properties of these processes are:

\begin{enumerate}
  \item With probability one,
  $Y_n$ assumes the values $y_k=y_{k,n}=k/(2^n-1)$ for $k=0,1,\ldots,2^n-1$, and
the difference between the realizations of $Y_n$ is a constant; more precisely,
 $y_{l+1}-y_l = \frac{l+1}{2^n-1}-\frac{l}{2^n-1} = 1/(2^n-1)$ for $l=0,\ldots,2^n-2$.

  \item ${\bf B}_n={\bf b}_n=(b_{n-1},\ldots,b_1,b_0)$ is a binary vector of length $n$. Thus $b_{j-1}$ may be interpreted as the $j$-th digit in the binary representation of a natural number $k \in \{0,\ldots,2^n-1\}$, i.e., $k=\sum_{j=0}^{n-1} b_j 2^j$, and the probability $p_k$ at the point $y_k$ is given by
  $$
  p_k = p^{\#1} (1-p)^{\#0} =  p^{\sum_{i=0}^{n-1} b_i}(1-p)^{n-\sum_{i=0}^{n-1} b_i} \ge 0 ,
  $$
where $\#1$ and $\#0$ denote the number of ones and zeros in ${\bf b}_n$, respectively. In particular, every $p_k$ can be written in the form $p_k=p^l (1-p)^{n-l}$ with some $l\in \{0,\ldots,n\}$.

  \item More explicitly, the distributions of ${\bf B}_n$, $E(S_n|{\bf B}_n)$, and ${Y}_n$ are
$$
\begin{array}{|r|lr|c|c|}\hline
E(S_n|{\bf b}_n)   & (k)_2 \Leftrightarrow & {\bf b}_n =(b_{n-1},\ldots,b_0)   & y_{k,n} & p_k \\
= (k)_{10}  & & & = k / (2^n-1) &  \\\hline
 0 & 0 & (0,\ldots,0) & 0 & (1-p)^n \\
 1 & 1 & (0,\ldots,0,1) & 1/(2^n-1) & p(1-p)^{n-1} \\
 2 & 10 & (0,\ldots,0,1,0) & 2/(2^n-1) & p(1-p)^{n-1} \\
 3 & 11 & (0,\ldots,0,1,1) & 3/(2^n-1) & p^2(1-p)^{n-2} \\
 4 & 100 & (0,\ldots,0,1,0,0) & 4/(2^n-1) & p(1-p)^{n-1} \\
 {\ldots}&{\ldots}&{\ldots}&{\ldots}&{\ldots} \\
 2^n-5 & 1\ldots 1011 & (1,\ldots,1,0,1,1) & (2^n-5)/(2^n-1) & p^{n-1}(1-p) \\
 2^n-4 & 1\ldots 100 & (1,\ldots,1,0,0) & (2^n-4)/(2^n-1) & p^{n-2}(1-p)^2 \\
 2^n-3 & 1\ldots 101 & (1,\ldots,1,0,1) & (2^n-3)/(2^n-1) &  p^{n-1}(1-p) \\
 2^n-2 & 1\ldots 10 & (1,\ldots,1,0) & (2^n-2)/(2^n-1) & p^{n-1}(1-p) \\
 2^n-1 & \underbrace{1\ldots 1} & (1,\ldots,1) & 1 & p^{n} \\
 & n \; \mbox{times} &&& \\\hline
\end{array}
$$
\item There is the representation
\begin{equation}\label{central}
  (2^n -1) y_{k,n} = E(S_n|{\bf b}_n)=\sum_{j=1}^{n} b_{j-1} 2^{j-1} \in \{0,\ldots,2^n-1\} .
\end{equation}

\end{enumerate}
  \end{thm}

Proof: (iv). Note that $X_i$ is neither Bernoulli nor integer-valued, and thus $S_n$ can be any real number. However, conditional on ${\bf B}_n={\bf b}_n$, with probability one, the expected sum $E(S_n|{\bf b}_n)$ is a natural number, since the typical contribution of an $X_i \sim H_c$ to the sum $S_n$ is either one if $c=1$ or zero if $c=0$. Since the $j$-th subsample consists of $2^{j-1}$ random variables having the same distribution, the right hand side of equation (\ref{central}) follows. The left hand side is due to definition.

Thus the random variable $E(S_n | {\bf B}_n) = \sum_{j=1}^{n} B_{j-1} 2^{j-1}$ takes values in $\{0,\ldots, 2^n -1\}$ which implies (i). (ii) is due to construction. This or the binomial theorem yield
$
\sum_{k=0}^{2^n-1} p_k = \sum_{l=0}^n \binom{n}{l} p^l (1-p)^{n-l} =1.
$ Finally, (i) and (ii) imply (iii). $\;\;\;\; \diamondsuit$

It is crucial that the populations differ in expectation ($\mu (H_0) \neq \mu (H_1)$), i.e., that $\bar{X}_n$ may fluctuate between two `centres of gravity'. The assumption $\mu(H_c)=c, (c=0,1)$ simplifies\label{simplify} the formal treatment considerably, since it leads to equation (\ref{central}) which states that $c=b_{j-1} \in \{0,1\}$ may be interpreted as the index of the selected population $H_c$, and as the expectation $EX_c$ where $X_c \sim H_c$. Otherwise, a linear transformation would have to map $\mu(H_0)$ to zero and $\mu(H_1)$ to one, which would add an additional - yet unnecessary - layer of complexity.

\vspace{2ex}
We say that $Y_n$ has a {\it weaver's distribution}, $Y_n \sim W(n,p)$, with parameters $n$ and $p$. Since powers of two play a major role, `binary distribution' would also be a suitable choice - much in line with `Bernoulli' and `binomial' distributions, which are closely related.

\vspace{2ex}
\begin{thm}\label{pascal}
(The `geometric triangle')

Given the assumptions and the notation of the last theorem, let ${\bf b}_n={\bf s}_{ij}$ be a vector with exactly $i$ ones and $j$ zeros, such that $i+j=n$. Moreover, set $f= p / (1-p)$.
\begin{enumerate}
  \item  The probabilities of the concatenated vectors $({\bf s}_{ij},1), (1,{\bf s}_{ij}), ({\bf s}_{ij},0)$, and $(0,{\bf s}_{ij})$ are: $$
     \frac{p({\bf s}_{ij},1) }{ p({\bf s}_{ij},0)}=\frac{p(1,{\bf s}_{ij})}{p(0,{\bf s}_{ij})} = \frac{ p^{i+1}(1-p)^j }{ p^i(1-p)^{j+1}} = \frac{p}{1-p} = f
    $$
    In particular, $ p_{l+1} / p_{l}=p/(1-p)=f$ for any
  two adjacent realizations $y_{2l}, y_{2l+1}$, where $l=0,1,\ldots,2^{n-1}-1$. The probabilities $p(\cdot)$ of the concatenated vectors $(0,1,{\bf s}_{ij}), (1,0,{\bf s}_{ij})$, etc., are
    $$
     \frac{p(0,1,{\bf s}_{ij})}{p(1,0,{\bf s}_{ij})} = \frac{p(0,{\bf s}_{ij},1)}{p(1,{\bf s}_{ij},0)} =\frac{p({\bf s}_{ij},0,1)}{p({\bf s}_{ij},1,0)} = \frac{ p^{i+1}(1-p)^{j+1} }{ p^{i+1}(1-p)^{j+1}} = 1
    $$
\item
For $n=1,2,\ldots$, $p_0 = p_0(n) =(1-p)^n$ is the probability that only $H_0$ is chosen, and $p_k = p_0 \cdot f^{\#1}$ for $k=0,\ldots,2^n-1$, where, again, $\#1$ is the number of ones in the binary representation of $k$. This means that the vector of probabilities ${\bf p}_n = (p_0,p_1,\ldots,p_{2^n -1})$ can be written as follows:
    \begin{eqnarray*}
    {\bf p}_n &=& p_0 \cdot (1;f;f,f^2;f,f^2,f^2,f^3;f,f^2,f^2,f^3,f^2,f^3,f^3,f^4;\ldots; \\
    && f,f^2,f^2,f^3,\ldots,f^{n-2},f^{n-1},f^{n-1},f^n) = p_0 \cdot {\bf f}_n
     \end{eqnarray*}

\item More explicitly, with ${\bf p}_0=1$, the vector ${\bf f}_n$ has dimension $2^n$ and obeys the recursive relation ${\bf f}_0=1$, and ${\bf f}_n = ({\bf f}_{n-1},f \cdot {\bf f}_{n-1})$ for $n=1,2,\ldots$ Thus its components can be calculated with the help of the following scheme, which may be interpreted as a geometric version of Pascal's triangle.\footnote{Pascal named his triangle ``triangle arithmetique.'' Thus, at least in French, it is straightforward to name the above multiplicative structure ``triangle geometrique.'' Since row $n$ has $2^n$ entries, the geometric triangle is a `real' triangle on the ld scale.}
    $$
    \begin{array}{lccccccccccccccc}
   n=0: &&&&& &&& 1  &  & &&&&& \\
   n=1: &&&& 1 &&&& | &&&& f &&&\\
   n=2: && 1 && | && f && || && f && | && f^2 & \\
   n=3: & 1 & | & f & || & f & | & f^2 & ||| & f & | & f^2 & || & f^2 & | & f^3 \\
   & \multicolumn{15}{c}{\ldots} \\
    \end{array}
    $$
    Every row has $2^n$ entries. Note that the left and the right of every $|$ are `separated' by the factor $f$ in the following sense: First $\left[ | \right]$, $1/f=f/f^2= f^2 / f^3 =\ldots$, or, equivalently, $1 \cdot f = f; f \cdot f = f^2; f^2 \cdot f = f^3$, etc. Second $\left[||\right]$, $(1,f) \cdot f = (f, f^2); (f,f^2) \cdot f = (f^2,f^3),(f^2,f^3)\cdot f = (f^3,f^4)$, etc. Third $\left[|||\right]$, $(1,f,f,f^2) \cdot f = (f,f^2,f^2,f^3);  (f,f^2,f^2,f^3) \cdot f = (f^2,f^3,f^3,f^4);$ etc.

    \item One may construct successive rows of (iii) in a rather elementary way: Start with a single 1 in the very first row. Then, fork every entry of row $n$ into two, by multiplying each entry with $1$ and $f$ upon moving left and right, respectively. It is quite remarkable that this local (cascade) view is equivalent to the global (weaving) view taken in the definition.\footnote{It may be noted that the `weaver' is similar to the `baker' in dynamic system theory. In particular, in both cases a locally defined transformation is closely related to global patterns. Theorem \ref{mandelbrot} connects the stochastic and the dynamic points of view explicitly.}

\item Applying the logarithm base $f$ to every entry of the geometric triangle yields the exponents, i.e., the following numbers:
    $$
    \begin{array}{l|cccccccccccccccccr}
  n && & &&&&  & &  &  &  & &&&&& & \mbox{Sum}\;\; s_n  \\\hline
   0 &&& &&&&  & &  & 0 &  & &&&&& & 0 \\
  1 & && &  & & 0 &&&& | &&&& 1 && & & 1\\
 2 &&&& 0 && | && 1 && || && 1 && | && 2 && 4\\
 3 & &&  0 & | & 1 & || & 1 & | & 2 & ||| & 1 & | & 2 & || & 2 & | & 3 & 12 \\
    \multicolumn{19}{c}{\ldots} \\
    \end{array}
    $$

In general, $s_0=0$, and $s_{n+1}=2 s_n +2^n$ for $n=0,1,\ldots$ That is, one obtains the sequence $0,1,4,12,32,80,192,448,1024,2304,\ldots$
\end{enumerate}
\end{thm}

Proof: (i) is proven in the statement of the theorem. However, (i) is also obvious, since the positions of the numbers 0 and 1 are irrelevant for the probabilities in question. In particular, for $k=0,2,\ldots,2^n-2$, the binary representations of $k$ and $k+1$ differ in exactly one position.

(ii) Using Theorem \ref{weaver} (ii), one obtains immediately

$p_k=p^{\#1} (1-p)^{\#0} = p^{\#1} (1-p)^{n-(\#1)} = (1-p)^n \frac{p^{\#1}}{(1-p)^{\#1}} = p_0 f^{\#1}$

(iii) is a consequence of self-similarity. Since the binary representations of 0 and $2^{n-1}$, and of 1 and $2^{n-1}+1$, etc., differ only by a single one,
\begin{eqnarray*}
  {\bf p}_{n} &=& (p_0,\ldots,p_{2^{n-1}-1};p_{2^{n-1}},\ldots,p_{2^{n}-1}) = (p_0,\ldots,p_{2^{n-1}-1};f p_0, f p_1 , \ldots,f p_{2^{n-1}-1}) \\
  &=&  ({\bf p}_{n-1}, f {\bf p}_{n-1})  = (p_0{\bf f}_{n-1}, f p_0 {\bf f}_{n-1} ) = p_0 ({\bf f}_{n-1}, f {\bf f}_{n-1} )
\end{eqnarray*}
Since, again by (ii), also ${\bf p}_{n} = p_0 {\bf f}_{n}$, the desired result follows.

One may also prove (iii) by induction on $n$: First, $p_1 = f p_0$, and thus $(p_0,p_1)=(p_0,f p_0) = p_0 (1,f)$.
 Second, the binary representation of any $k \in \{0,\ldots,2^n-1\}$ is a vector ${\bf b}_n = (b_{n-1},\ldots,b_0)$. Let $\#1$ be the number of ones in ${\bf b}_n$. With probability $1-p$, the next selection leads to $(0,{\bf b}_n)$, and with probability $p$ this selection results in $(1,{\bf b}_n)$. Since in the first case, the number of ones does not change, and in the second case, the number of ones increases by one, we obtain on the one hand (to the left), $p_{i,n+1}  = p_{0,n+1} f^{\#1} = (1-p)^{n+1} f^{\#1} = (1-p) p_{0,n} f^{\#1} =(1-p) p_{i,n}$ for $0 \le i \le 2^n-1$. This is tantamount to ${\bf f}_n$ being reproduced as the first half of ${\bf f}_{n+1}$. (Upon moving from $n$ to $n+1$, the exponent of $f$ does not change.) On the other hand (to the right), we obtain $p_{i,n+1}  = p_{0,n+1} f^{(\#1) +1} =  (1-p)^{n+1} f^{\#1} p/(1-p)= p(1-p)^n f^{\#1} =p p_{0,n} f^{\#1} = p p_{i,n}$ for $2^n \le i \le 2^{n+1}-1$. The additional factor $f$ means that the second half of ${\bf f}_{n+1}$ has to be $f \cdot{\bf f}_n$.

(iv) The proof is by induction on $n$. For $n=0$ there is nothing to prove, and the equivalence is obvious for $n=1$. By the inductive assumption, the vector occurring on line $n$, having length $2^n$, has the form ${\bf w}_n =({\bf l}_{n-1},{\bf r}_{n-1})=({\bf l}_{n-1}, f \cdot {\bf l}_{n-1})$ where ${\bf l}_{n-1}$ is a vector of length $2^{n-1}$. In other words, $r_k/l_k = f$ for $k=1,\ldots,2^{n-1}$.

    Local splits (see the definition given in the statement of the theorem) produce a vector ${\bf w}_{n+1}$ of length $2^{n+1}$. 
    Since, locally, a step to the left reproduces the numbers, and a step to the right multiplies any two entries on tier $n$ with the same factor $f$, we also have, because of the inductive assumption, $w_{2^n+k}/w_{k}=f$ for $k=1,\ldots,2^{n}$. 
    Therefore ${\bf w}_{n+1} = ({\bf l}_n,f \cdot {\bf l}_n)$.

(v) Straightforward induction on $n$ yields the recursive formula. $\diamondsuit$

\vspace{1ex}
Note that the multiplicative triangle lies at the heart of the observation that ``the best known multifractal constructions use multiplicative operations'' (\citet{ma99}, p. 32).

\vspace{3ex}
\begin{thm}\label{weaver2}

(Further properties of the weaver's distribution).

Given the assumptions and the notation of Theorem \ref{weaver}, one obtains
\begin{enumerate}
    \item The probabilities corresponding to row $n$ can be constructed by the following simple scheme:
     $$
     \hspace{-15ex}
    \begin{array}{ccccccccccccccc}
    &&&&  & &  & 1 &  & &&&&& \\
    &&& 1-p &&&& | &&&& p & &&\\
    & \text{\begin{small}$(1-p)^2$\end{small}} && | && \text{\begin{small}$p(1-p)$\end{small}}  && || && \text{\begin{small}$p(1-p)$\end{small}}  && | && \text{\begin{small}$p^2$\end{small}}  & \\
\text{\begin{tiny}$(1-p)^3$\end{tiny}}   & | & \text{\begin{tiny}$p(1-p)^2$\end{tiny}}   & || & \text{\begin{tiny}$p(1-p)^2$\end{tiny}}
 & | & \text{\begin{tiny}$p^2(1-p)$\end{tiny}}  & ||| & \text{\begin{tiny}$p(1-p)^2$\end{tiny}} & | & \text{\begin{tiny}$p^2(1-p)$\end{tiny}}
 & || & \text{\begin{tiny}$p^2(1-p)$\end{tiny}}  & | & \text{\begin{tiny}$p^3$\end{tiny}}  \\
    \multicolumn{15}{c}{\ldots} \\
    \end{array}
    $$
 Global interpretation [weaving]: ${\bf p}_{n+1} = ((1-p){\bf p}_n,p {\bf p}_n)$. Local interpretation [Bernoulli cascade]: Start with mass 1 in the very first (the zeroth) row. Then, fork every probability of row $n$ into two, by multiplying each entry with $1-p$ (on the left) and $p$ (on the right), respectively.

    \item The sequence $p_0,fp_0,f^2p_0,\ldots$ increases (decreases), the mode occurs at one (zero), and the median is larger (smaller) than $1/2$ if $p>1/2$ $(p<1/2)$. If $p=1/2$, all probabilities coincide, i.e. we obtain the discrete uniform distribution on the values $y_k=k/(2^n-1)$; $p_k=1/2^n$ for $k=0,1,\ldots,2^n-1$.
\item Symmetry: Suppose $Y \sim W(n,p)$ and $Y^{'} \sim W(n,1-p)$. Then $P(Y=y_k) = P(Y^{'}=y_{2^n-1-k})$ for $k=0,\ldots,2^n -1$. 

    \item Distribution function $F$ of $W(n,p)$: For all $n\ge 0$ and $k=0,\ldots,2^n$ define $v_{k,n}=k/2^n$. For every fixed $n$, the mass left and right of $v_{k,n}$ $(0<k<2^n)$ is constant for every $m \ge n$, and so is the value of $F(v_{k,n})$. In particular, $F(v_{1,1})=F(1/2)=(1-p)$ for all $n \ge 1$; $F(v_{1,2})=F(1/4)=(1-p)^2$, $F(v_{3,2})=F(3/4)=1-p^2$ for all $n\ge 2$; $F(v_{1,3})=F(1/8)=(1-p)^3$; $F(v_{3,3})=F(3/8)=(1-p)^2+p(1-p)^2$, $F(v_{5,3})=F(5/8)=(1-p)+p(1-p)^2$, $F(v_{7,3})=F(7/8)=1-p^3$ for all $n\ge 3$, etc.
        \item The total mass $p_k$ in every interval $[v_{k,n},v_{k+1,n}]$ ($k=0,\ldots,2^n-1$) remains the same for all $m \ge n$. For $m=n$ it is located at the point $y_k=y_{k,n}=k/(2^n-1)$. In the interest of consistency let $y_{0,0}=p$ and $p_0=1$ if $n=0$.

            Thus $W(n,p)$ may be interpreted as a discretisation of the density in the corresponding classical $p$ model.
            \item Distribution of the jumps (stick heights): $F_n$ has $2^n$ points of discontinuity. If $p=1/2$ there is a constant jump height $h=1/2^n$. Otherwise, there are $n+1$ different jump sizes, given by $h_j=p^j(1-p)^{n-j}$ for $j=0,\ldots,n$, having a binomial distribution. That is, there is 1 jump of size $h_0=(1-p)^n$, there are $\binom{n}{1}=n$ jumps of size $h_1=(1-p)^{n-1}p$, etc.
      \end{enumerate}
\end{thm}

Proof: \begin{enumerate}
         \item For $n=1,2,\ldots$, we have $p_0=p_0(n)=(1-p)^n$ for the leftmost probability (only $H_0$ is selected). Applying the geometric triangle yields the result.
         \item We have $p<1/2 \Rightarrow f>1$. Thus the mass in $y_1$ exceeds the mass in $y_0=0$ by the factor $f$, and the result follows straightforwardly.
\item Exchanging the roles of zeros and ones, and replacing $p$ by $1-p$ yields the same distribution. In other words: The reflection of $W(n,p)$ across the axis of symmetry $y=1/2$ is $W(n,1-p)$.
    \item follows immediately from the geometric triangle. Geometrically speaking, the unit interval on the horizontal axis is successively halved. At the same time, the unit interval on the vertical axis is successively divided according to the ratio $f$. Thus, for finite $n \ge 1$, one obtains a step function with $2^n$ jumps.
\item holds because of the local interpretation of the geometric triangle: Each split can be interpreted as distributing the mass $p_k$ in $y_k$ to the points $y_{2k,n+1}$ and $y_{2k+1,n+1}$ in that same interval. Graphically, the stick of height $p_k$ in $y_{k,n}$ is broken into two sticks of heights $(1-p)p_k$ and $p \cdot p_k$, located in $y_{2k,n+1}$ and $y_{2k+1,n+1}$, respectively.
    \item is due to construction.
$\diamondsuit$
  \end{enumerate}


Note that there are two kinds of scale: the first could be named `discrete time,' i.e., the total number of observations $2^n$, the second would be `logarithmic time,' that is, the number of selections, $\ld 2^n =n$.

It is also instructive to compare the above structure to the received $\alpha$ model (see, in particular, \citet{lo13}, pp. 65-70, with the number of descendants being $\lambda_0$). If $\lambda_0 =2$, both models are built on a binary cascade, i.e., after $n$ steps, there are $2^n$ leaves. Like the $p$ model, the $\alpha$ model distributes the mass assigned to some interval $I_k=(k/2^n,(k+1)/2^n)$, $k=0,\ldots,2^n-1$ uniformly on that interval. Moreover, the initial mass in some parent vertex is split according to a binary random variable $B$ that assumes the value $\gamma_+ >0$ (a `boost') with probability $p=2^{-c}$ and the value $\gamma_{-}<0$ (a `decrease') with probability $1-p=1-2^{-c}$, where $c$ is a positive parameter. In particular, if the mass at the root is one,
$$
\begin{array}{rcccl}
 && 1 &&  \\
(1-p)=1-2^{-c} &\swarrow &  & \searrow & p = 2^{-c} \\
 2^{\gamma_{-}} &&&&  2^{\gamma_{+}} \\
\end{array}
$$
and one assumes $EB=2^{\gamma_{+}}\cdot 2^{-c} + 2^{\gamma_{-}} (1-2^{-c})=1$ in order to get a canonical cascade.
The illustrations ibd., p. 68, show that this model produces a `randomized $p$ density', i.e., a function that is equal to a constant $c_k >0$ on every interval $(k/2^n,(k+1)/2^n), k=0,\ldots,2^n-1.$ However, the values $c_k$ are random, since
$n$ steps amount to exactly $\sum_{j=1}^{n} 2^{j-1} = 2^n-1$ binary decisions, and each decision means to split the mass $m$ found at some parent vertex in the way just described, i.e., with probability $p$ a descendant inherits mass $m \cdot 2^{\gamma_+}$, with probability $1-p$ it receives mass $m \cdot 2^{\gamma_{-}}$. Thus one obtains a particular realization of a sequence / cascade of randomly occurring boosts and decreases.

Weaving is different: First, we are interested in all possible realizations. Therefore, after $n$ steps, the total mass of 1 has been distributed among the $2^n$ leafs $y_{k,n}$ of a binary tree. So, secondly, there is no density but a discrete distribution. Third, and perhaps most importantly, there have been $n$ `major' decisions that chose a parent distribution ($H_0$ or $H_1$) for each of the stages $j=1,\ldots,n$. Based on the latter distributions, there have also been $2^n-1$ `minor' decisions. That is, given $j$ (and thus either $H_0$ or $H_1$), one has obtained a group of $2^{j-1}$ independent realizations $x_{2^{j-1}},\ldots,x_{2^{j}-1}$ from the parent distribution in charge. In other words, the cascade corresponding to weaving is hierarchical: $n$ binary decisions select the composition of the sample $X_1,\ldots,X_{2^n-1}$ (exactly one out of $2^n$ possible compositions, since $j=1,\ldots,n$). Given this, in each group $j$, one considers the sum $T_j= \sum_{i=2^{j-1}}^{2^j -1} X_i$ and its realization $t_j= \sum_{i=2^{j-1}}^{2^j -1} x_i$.

A different interpretation would be that the basic building block of an (unstandardized) Bernoulli cascade is
$$
\begin{array}{rcccl}
 && S_{j-1} &&  \\
(1-p) &\swarrow &  & \searrow & p \\
X_i \sim H_0  &&&& X_i \sim H_1 \\
S_j = S_{j-1} + T_j  &&&& S_j = S_{j-1} + T_j = S_{j-1} + \sum_{i=2^{j-1}}^{2^j -1} X_i \\
\end{array}
$$
where $j=1,2,\ldots$ and $S_0=0$. Notice that the crucial feature of this building block and the corresponding cascade is a particular kind of summation.

\section{Moments}

\begin{thm}\label{expectedvalue}

{\it (Expected value)}.
Let $Y_n \sim W(n,p)$. Then, for every $n \ge 1$, $EY_n = p$.
\end{thm}

Proof: Equation (\ref{central}) implies
$$
     E(Y_n) = E \left( \sum_{j=1}^{n} B_{j-1} 2^{j-1})/(2^n-1) \right)
        = \frac{1}{(2^n-1)} \sum_{j=1}^{n} 2^{j-1}  E( B_{j-1}) = p ,
$$
since $\sum_{j=1}^{n} 2^{j-1} = 2^n -1$. $ \diamondsuit$


After the first step, the distribution of the conditional expected values is $B(p)$. For any random variable $X$ with values in the unit interval, and $EX=p$, this distribution has maximum variance $p(1-p)$. Upon weaving, probability mass is successively concentrated within the unit interval, and thus variance decreases. On the other hand, every bifurcation may increase the variance term.

Both effects combined result in a (net) monotone decrease of variance up to a certain point. Moreover, there is a limit variance $\sigma^2= c p (1-p)$ with $0<c <1$.

\begin{thm}\label{variancetheorem}
{\it (Variance).}
Let $Y_n \sim W(n,p)$. Then the variance of this r.v. is
\begin{equation}\label{variance}\sigma^2(Y_n) =\frac{\sum_{i=0}^{n-1} 2^{2i}}{(2^n-1)^2} p (1-p) = \frac{4^n -1}{3(2^n -1)^2} p (1-p) \end{equation}
\end{thm}

Proof: Representation (\ref{central}) and the independence of the $B_j$ imply \begin{eqnarray*}
       \sigma^2(Y_n) &=& \sigma^2(\sum_{j=1}^{n} B_{j-1} 2^{j-1})/(2^n-1)) = \frac{1}{(2^n-1)^2} \sum_{j=1}^{n} \sigma^2( B_{j-1} 2^{j-1})  \\
        &=& \frac{1}{(2^n-1)^2} \sum_{j=1}^{n} (2^{j-1})^2  \sigma^2( B_{j-1}) = \frac{p (1-p)}{(2^n-1)^2} \sum_{j=1}^{n} 2^{2(j-1)} ,
     \end{eqnarray*}
and $\sum_{i=0}^{n-1}2^{2i} = \sum_{i=0}^{n-1} 4^{i} = \frac{1-4^n}{1-4} = \frac{4^n -1}{3} . \;\;\;\;\diamondsuit$

Note that in the quaternary system $((4^n-1)/3)_4 = (1\ldots 1)_4$ with $n$ ones in the last bracket.
Moreover, the latter sums have a nice interpretation:
$$
\begin{array}{c|ccccc|r}
 {\bf M}  & 1 & 2 & 4 & \ldots & 2^{n-1} & \sum \\\hline
 1 & {\bf 1} & 2 & 4 & \ldots & 2^{n-1} & 2^n-1  \\
  2 & 2 & {\bf 4} & 8 & \ldots & 2^n & 2(2^n-1) \\
  4 & 4 & 8 & {\bf  16} & \ldots & 2^{n+1} & 4(2^n-1) \\
  \ldots & \ldots & \ldots & \ldots & \ldots & \ldots & \ldots \\
  2^{n-1} & 2^{n-1} & 2^n & 2^{n+1} & \ldots & {\bf 2^{2(n-1)}} & (2^{n-1})(2^n-1)\\\hline
 \sum & 2^n-1  & 2(2^n-1) & 4(2^n-1) & \ldots & (2^{n-1})(2^n-1) & (2^n-1)^2
\end{array}
$$

\vspace{1ex}
That is, $tr({\bf M})=\sum_{j=1}^{n} 2^{2(j-1)}$ is the trace of the matrix ${\bf M}$, $(2^n-1)^2 - tr({\bf M})= \sum_{j=0}^{n-1} 2^j (2^n-1-2^j) = \frac{2}{3} (2-3 \cdot 2^n + 4^n)$ is the sum of all elements not in the main diagonal of ${\bf M}$, and $\sum_{j=1}^{n} 2^{2(j-1)}/(2^n-1)^2$ is the proportion of the trace relative to the total sum $(2^n-1)^2$.

Also note that the numerator shows an additive analogue to factorials: For factorials, $n! = (n-1)! \cdot n$ holds. For the numerator, we have $\sum_{i=0}^{n} (2^{i})^2 = \sum_{i=0}^{n-1} 2^{2i} + 2^{2n}$.

\begin{cor}\label{moments}
$EY_n^2$ exists, and so do all higher moments $EY_n^j$ for $j \ge 1$.
\end{cor}

Proof: For fixed $n$, all realizations $y_k$ are in the unit interval. Thus $y_k \ge y_k^2 \ge y_k^3 \ge \dots$, with strict inequality if $0 < y_k < 1$. Therefore $0 < EY_n^i < EY_n^j$ if $i > j$. $\diamondsuit$

\begin{thm}\label{extension}
  Suppose $B_{k,n} \sim B(y_{k,n})$ where $y_{k,n}$ are the realizations of $Y_n \sim W(n,p)$. Denote by $Z_n = B_{k,n} \circ Y_n$ the random variable defined by replacing the realization $y_{k,n}=k/(2^n-1)$ with the Bernoulli r.v. $B_{k,n}$ $(k=0,\ldots,2^n -1)$.

  Then $Z_n$ is $\{0,1\}$ - valued, $EZ_n =p$ and $\sigma^2(Z_n)= p (1-p)$.
\end{thm}

Proof: Since each $y_{k,n}$ is mapped to 1 w.p. $y_{k,n}$ and to 0 w.p. $1-y_{k,n}$, $Z_n =B_{k,n} \circ Y_n$ assumes the values zero and one. Moreover, since $EB_{k,n}=y_{k,n}$ the location of $Y_n$'s distribution is the same as that of $Z_n$. Thus $Z_n \sim B(p)$ and the variance follows. $\diamondsuit$

For fixed $n$, one might think of the collection of all $B_{k,n}$ $(k=0,\ldots,2^n -1)$ as a family of ``dual distributions'' to $W(k,n)$, mapping that distribution to the set $\{ 0,1 \}$ without changing the center of gravity.
Since the point $y_{k,n}$ and its mass are distributed to two values, `splitting' might be an appropriate term for this operation. However, since the crucial idea is that all $y_{k,n}$ map to the set $\{0,1\}$, `merging' is even more appropriate.

A different interpretation would be that upon constructing the Binomial, splitting and merging occur together (cannot be separated) when moving from $n$ to $n+1$. The Bernoulli cascade also starts with unit mass at point $p$. However, after $n$ bifurcations, it creates a $W(n,p)$ that distributes mass among $2^n$ points. Hereafter, `all the merging' occurs in a single step, i.e., merging is postponed until the end.

Due to
Theorems \ref{variancetheorem} and \ref{extension}, one may decompose the total variance among weaving and merging:
\begin{eqnarray*}
  \sigma^2(Z_n)  &=& p(1-p) = \frac{tr({\bf M})}{(2^n -1)^2} p (1-p)  + \frac{(2^n -1)^2 -tr({\bf M}) }{ (2^n -1)^2} p (1-p)  \\
   &=& \frac{\sum_{i=0}^{n-1} 2^{2i}}{(2^n-1)^2} p (1-p) - \frac{\sum_{j=0}^{n-1} 2^{2j}  }{(2^n-1)^2} p (1-p) + \frac{\sum_{j=0}^{n-1} 2^j (2^n-1) }{(2^n-1)^2} p (1-p) \\
   &=& \sigma^2(Y_n) + \frac{\sum_{j=0}^{n-1} 2^j (2^n-1-2^j) }{(2^n-1)^2} p (1-p) \\
    &=& \sigma^2(E(Z_n | {\bf B}_n)) + E (\sigma^2 (Z_n | {\bf B}_n))
\end{eqnarray*}

Note that the last term in the second line is just the last column of matrix ${\bf M}$. More importantly, in
the last line, we interpret $Z_n$ as a particular mixture of $H_0$ and $H_1$, whose total variance may be decomposed into the ``variance of the conditional expected values'' (the first term, i.e., the variability in the $y_{k,n}$) and the ``expected conditional variance'' (the second term, i.e., the mean variance of the $B_{k,n}$).

Some concrete values may be helpful:
$$\label{concrete}
\begin{array}{|l|r|r|r|r|c|}\hline
&& Denom. & Weaving & Merging & Proportions \\\hline
  n & 2^n-1 & (2^n-1)^2 & (4^n-1)/3 & 2(4^n -3 \cdot 2^{n}+ 2)/3 &  \\\hline
 1 & 1 & 1 & 1 & 0 & 1 \leftrightarrow 0 \\
 2 & 3 & 9 & 5 & 4 & \frac{5}{9}=0.{\bar 5} \leftrightarrow 0.{\bar 4} = \frac{4}{9} \\
 3 & 7 & 49 & 21 & 28 & 0.43 \leftrightarrow 0.57 \\
 4  & 15 & 225 & 85 & 140 & 0.3{\bar 7} \leftrightarrow 0.6{\bar 2}\\
 5 & 31 & 931 & 341 & 620 & 0.35 \leftrightarrow 0.65 \\
 6 & 63 & 3969 & 1365 & 2604 & 0.34 \leftrightarrow 0.66 \\
  \hline
\end{array}
$$

\section{Limit distribution}

\begin{lem}\label{limitvariance}
The limit of the variance term is $\frac{1}{3}p(1-p)$
\end{lem}

Proof: Considered as a function of $n$, $\sigma^2(Y_n)$ is monotonically decreasing. Since it is also nonnegative, it is clearly convergent. Thus we get
$$
\frac{\sigma^2(Y_n)}{p (1-p)} = \frac{\sum_{i=0}^{n-1}2^{2i} }{ (2^n -1)^2}= \frac{4^{n}-1}{3 (2^{2n}-2^{n+1}+1)}
   = \frac{(2^{2n}-1) / 2^{2n}}{3 (2^{2n}-2^{n+1}+1) / 2^{2n}} ,  \\
$$
which converges to $1/3$ if $n \rightarrow \infty$. $\diamondsuit$

   Since, due to Theorem \ref{weaver2}, the distribution function $F_n$ is well-known for all values $v(k,n)$, it is easy to pass to the limit. The limit function $F$ obviously is a distribution function.

\begin{thm}\label{limitdistribution} (The weaver's hem)

Let $Y$ be the limit of $(Y_n)$, defined by its distribution function $F=\lim_{n \rightarrow \infty} F_n$. For obvious reasons, the corresponding random variable and distribution $Y \sim W(p)$, might be named the {\it weaver's hem.}

Then $F$ is continuous, $EY=p$ and $\sigma^2(Y)=p(1-p)/3$. Except for the case $p=1/2$, when $W(n,1/2)$ is a discrete uniform distribution and thus $W(1/2) \stackrel{d}{=} U(0,1)$, (implying that $F(x)=x$ for $0<x<1$), $F$ has no density with respect to Lebesgue measure.
\end{thm}

Proof: Using the notation of Theorem \ref{weaver2},
for fixed $n$, all mass is concentrated at the points $y_{k,n}=k/(2^n -1)$, $(k=0,\ldots,2^n-1)$, and the jump heights there (Theorem \ref{weaver2} (vii)) go to zero if $n \rightarrow \infty$. Thus $F$ has to be continuous.

Because of $EX=\int_0^1 (1-G(x)) dx$ for any distribution function $G$ on the unit interval, and $F_n\rightarrow F$, we also have $EY=p$ for the weaver's hem. An analogous argument for the second moment and Lemma \ref{limitvariance} yield $\sigma^2(Y)=p(1-p)/3.$

   Rather heuristically, if $p > 1/2$, consider the interval $[0,1/2[$. The mass of $1-p$ available there is shifted to the left. Thus the distribution function grows rapidly at first, but hardly grows near $1/2$. Now consider the interval $]1/2,1]$. Because a mass of $p$ is available there and systematically shifted to the left, the distribution function grows rapidly near $1/2$, but very slowly near $1$. Thus the distribution function has a sharp point at $1/2$ and cannot be differentiated there. The same holds for all $v(k,n)$. Since the set of these points lies dense in the unit interval, there should be no density.

   Formally, consider the interval $[v_{k,n},v_{k+1,n}]$ about $y_k=y_{k,n}$. For fixed $n$, this interval has length $v_{k+1,n}-v_{k,n} = (k+1-k)/ 2^n = 1/2^n$.
   By Theorem \ref{pascal} (ii), the density in the neighbourhood of $y_k$ is given by
   \begin{equation}\label{dichte}
   g_{k,n}= 2^n p_k = 2^n p_0(n) f^{\#1} = 2^n (1-p)^n f^{\#1} = 2^n  p^{\#1} (1-p)^{\#0},
   \end{equation}
   where $\#0$ and $\#1$ are the number of zeros and ones in the binary representation of $k$, respectively.
   If $p=1/2$, $g_{k,n}=1$, and thus $W(1/2)$ is the uniform distribution on $[0,1]$.

   In general, compare Equation (\ref{dichte}) and the classical De Moivre-Laplace theorem. In the latter case, one considers $\binom{n}{k}p^k (1-p)^{n-k}$, which approaches a limit $b \in (0, \infty)$, since the convergence of $p^k(1-p)^{n-k}$ toward zero is counterbalanced by a sequence that goes to infinity at the same speed, i.e., an appropriate binomial coefficient (also depending on $n$ and $k$).

   Here, every iteration ($n \rightarrow n+1$) doubles the number of values $y_k$, and thus the first factor is $2^n$ instead of $\binom{n}{k}$. Moreover, due to Theorem \ref{weaver2}, every $y_{k,n}$ is the starting point of a cascade, i.e., a sequence of local bifurcations in the corresponding interval $[v_{k,n} ; v_{k+1,n}]$. After one iteration, the probabilities at $y_{2k,n+1}$ and $y_{2k+1,n+1}$, i.e., $(1-p)p_k$ and $p \cdot p_k$, respectively, differ by the factor $f$. After $l$ iterations, the probabilities at the leftmost value $y_{2^l k,n+l}$ and the rightmost value $y_{2^l k + (2^l -1),n+l}$ differ by $f^l$. If w.l.o.g. mass is systematically shifted to the right $(p>1/2)$, we have $f>1$, and thus the ratio of these probabilities soon exceeds any bound. Even more so, $2^l (1-p)^l p_k \rightarrow 0$ and $2^l p^l p_k \rightarrow \infty$ in every interval $[v_{k,n};v_{k+1,n}]$ if $l \rightarrow \infty$. 
   Thus, there cannot be a limit density. $\diamondsuit$

\vspace{2ex}
Note that the `roughness' of the density (measured by $f^l$) grows at the same rate as the number of intervals. Thus $\ln f^n / \ln 2^n = \ln f/ \ln 2$ is a constant, the fractal dimension.



\vspace{1ex}
\begin{thm}\label{mandelbrot}
   The weaver's hem and Mandelbrot's `binomial measure' are equivalent.
\end{thm}

Proof: Mandelbrot's `binomial measure' is the limit of the $p$ model, splitting the mass (locally) according to the geometric triangle.
Thus, the $p$ model's Bernoulli cascade and weaving (see Theorem \ref{weaver2}, (vi)) assign the same mass to every interval $[v_{k,n} ; v_{k+1,n}]$. Since these intervals shrink to zero, the limit distributions have to coincide.
$\diamondsuit$

\vspace{2ex}
Of course, the last theorem and \citet{sa43} also imply that $W(p)$ has no density if $p \neq 1/2$.

The last theorem could be an example of a more general `sandwich principle'. That is, $W(p)$ is defined on a countable set and is the limit of a Bernoulli cascade. The cascade either starts with unit mass at the point $p$ and distributes that mass to an increasing number of points (this article); or the cascade starts with the Uniform on the unit interval, and distributes that mass to an increasing number of (shrinking) intervals (the received $p$ model). If both processes have the same `inheritance' preference $f=p/(1-p)$, they determine the same multifractal in the limit.

\section{The complete process}

So far, we have mainly considered the distribution of the (conditional) expected values, $Y_n = E({\bar X}_n | {\bf B}_n)$, or, equivalently, the case of two one-point distributions located in $\mu(H_0)$ and $\mu(H_1)$, respectively. Looking at ${\bar X}_n$, however, there is not just variance between the populations $H_0$ and $H_1$, but also within each of these populations, $\sigma^2(H_0) = \sigma_0^2$ and $\sigma^2(H_1)=\sigma_1^2$, say, contributing to the total variance.

For the unconditional process we obtain:

\begin{thm}\label{amproperties-neu} (Expected value and variance)

  With the assumptions of Theorem \ref{weaver} , $E{\bar X_n}  =p$
  and
  \begin{equation}\label{varianz}
  \sigma^2({\bar X}_n) = \sigma^2 (Y_n)+ \frac{p \sigma_1^2 + (1-p) \sigma_0^2}{2^n-1}
 \end{equation}
\end{thm}

Proof:

Since the random variables $X_i$ are organized in groups (subsamples), it is helpful to consider $T_j=\sum_{i=2^{j-1}}^{2^j -1} X_i$ $(j=1,\ldots,n)$ and the decomposition $S_n = \sum_{j=1}^{n} T_j = \sum_{j=1}^{n} \sum_{i=2^{j-1}}^{2^j -1} X_i$.

For the expected value observe that $E\bar{X}_n = ES_n / (2^n-1)$, and
$$
ES_n = E \left( \sum_{j=1}^{n} \sum_{i=2^{j-1}}^{2^j -1} X_i \right) = \sum_{j=1}^{n} \sum_{i=2^{j-1}}^{2^j -1} EX_i = \sum_{j=1}^{n} \sum_{i=2^{j-1}}^{2^j -1} p = (2^n-1) p ,
$$
since $EX_i = p \cdot \mu(H_1) + (1-p) \cdot \mu(H_0) =p$. (All random variables of some group $j$ have distribution $H_1$ with probability $p$, and $H_0$ with probability $1-p$. Thus these are also the corresponding probabilities of any single $X_i$ $(i=1,\ldots,2^n-1)$.)

$T_j$ is the sum of $2^{j-1}$ random variables, and selection $B_j$ chooses $H_1$ with probability $p$, and $H_0$ with probability $1-p$. Given $H_c$ $(c=0,1)$, the conditional variance of $T_j$ is $\sigma^2(T_j | H_c)=2^{j-1} \sigma^2(X_i|H_c)=2^{j-1} \sigma_c^2$, since $X_{2^{j-1}},\ldots,X_{2^j -1}$ are iid.

Moreover, $E(T_j|H_c)=2^{j-1} c$, and $ET_j = 2^{j-1}p$. A variance decomposition of $T_j$ yields
\begin{eqnarray*}
  \sigma^2(T_j) &=& p (E(T_j|H_1)-ET_j)^2 + (1-p) (E(T_j|H_0)-ET_j)^2  \\
  && + p \sigma^2(T_j|H_1) + (1-p) \sigma^2(T_j|H_0) \\
  &=& p (2^{j-1}-2^{j-1}p)^2 + (1-p) (0-2^{j-1}p)^2 + p 2^{j-1} \sigma_1^2 + (1-p) 2^{j-1} \sigma_0^2 \\
  &=& 2^{2(j-1)} p (1-p) (1-p + p) + 2^{j-1} p \sigma_1^2 +  2^{j-1} (1-p) \sigma_0^2
\end{eqnarray*}
which implies
\begin{eqnarray}\label{basic}
\nonumber
  \sigma^2 (\bar{X}_n ) &=& \sigma^2 \left( \frac{S_n}{2^n -1} \right) = \frac{1}{(2^n -1)^2} \sigma^2 (S_n) = \frac{1}{(2^n -1)^2} \sum_{j=1}^{n} \sigma^2(T_j)   \\\nonumber
  &=&  \frac{\sum_{j=1}^{n} 2^{2(j-1)}}{(2^n-1)^2} p(1-p) + \frac{ \sum_{j=1}^{n} 2^{j-1}}{(2^n-1)^2} p \sigma_1^2   + \frac{\sum_{j=1}^{n}2^{j-1}}{(2^n -1)^2} (1-p) \sigma_0^2 \\
  &=& \sigma^2(Y_n) + \frac{p \sigma_1^2 + (1-p) \sigma_0^2}{2^n -1} ,
\end{eqnarray}
completing the proof. $\diamondsuit$

\section{Asymptotic structures}

Consistently, the variance within the populations washes out quickly, and one obtains $\lim_{n \rightarrow \infty} \sigma^2(\bar{X}_n) = \sigma^2(Y)= p (1-p ) /3$. In other words, in the limit only the variance caused by exponential sampling (the selection process) is relevant, and within the populations, the LLN applies.
This means:

\begin{thm}\label{amlimit} (Limit distribution)
Given the assumptions and the notation of Theorem \ref{weaver}, if $H_0$ and $H_1$ both have finite variances, $W(p)$ is the limit distribution of the inhomogeneous (unconditional) stochastic process $({\bar X}_n)$.
\end{thm}

A major reason for the universality of the CLT is that particular properties of the random variables involved do not matter in the limit. That is, suitably standardized iid random variables $X_i$ always converge in distribution toward the Normal, if they have finite variances. Quite similar here: As long as the two populations have finite variances, exponential sampling (plus a suitable standardization) yields $W(p)$ in the limit. Again, specific features of the distributions in question - details - fade, and only the sampling process (i.e., the way, sampling is done) determines the final result.

One could also argue that the CLT holds, since Lindeberg's condition implies that the contribution of a single random variable to the sum $\sum_{i=1}^{n} X_i$ becomes arbitrarily small. (Similarly, common wisdom has it that the Normal occurs if {\it many} small, independent errors add up.) Contrary to this, it is quite characteristic for stable (limit) distributions that a single variable has a dominant influence on the whole (see, in particular, \citet{fe71}, pp. 172, 465, \citet{lo19}, p. 45, and the remarks of \citet{ma97, ma99} on ``wild chance''). The Weaver's hem is observed if there are two `forces' ($H_0$, $H_1$), but neither of them dominates. Therefore, both contribute to the asymptotic outcome, which is the Uniform if there is no preference $(p=1/2)$. Otherwise, since threads do not merge, the `compromise' is a multifractal structure.

Finally, another straightforward interpretation would be that Pascal's combinatorial pattern yields the Binomial for finite $n$, and the Normal in the limit. Analogously, the `geometric triangle' or the corresponding Bernoulli cascade yield $W(n,p)$ and $W(p)$.

Extending $W(n,p)$ the way we did in Theorem \ref{extension}, asymptotically leads to $Y \sim W(p)$ and a countable family of dual distributions $B_k$. Since $\sigma^2(Y)=p(1-p)/3$, one third of the total asymptotic variance between the populations is due to weaving, and the rest stems from forcing all threads to terminate in 0 or 1.

Of course, if the populations $H_0$ and $H_1$ are not too complicated, it is possible to study the process ${\bar X}_n$ in much more detail. For various extensions see the last section of \citet{sa19}.

\end{doublespace}

\end{document}